\documentclass[reqno,12pt]{amsart}
\usepackage{amsmath,amsthm,amssymb,amsfonts,amscd}
\usepackage{xypic}
\theoremstyle{plain} 
	\newtheorem{thm}{Theorem}[section]
	\newtheorem*{thm*}{Theorem}
	\newtheorem{cor}[thm]{Corollary}
	\newtheorem{lem}[thm]{Lemma}
	
	\newtheorem{prop}[thm]{Proposition}
	\newtheorem{conj}[thm]{Conjecture}
	\newtheorem*{conj*}{Conjecture}
	
\theoremstyle{definition}
	\newtheorem{defn}[thm]{Definition}

\theoremstyle{remark}
	\newtheorem{rem}[thm]{Remark}
	
	\newtheorem*{pf}{Proof}
\numberwithin{equation}{section}
\def\CC{{\mathbb C}}

\def\PP{{\mathbb P}}
\def\QQ{{\mathbb Q}}

\def\ZZ{{\mathbb Z}}

\def\D{{\mathcal D}}
\def\E{{\mathcal E}}

\def\O{{\mathcal O}}

\def\S{{\mathcal S}}
\def\T{{\mathcal T}}

\def \mf#1#2#3#4{
\xymatrix{
{#1}\  \ar@<0.4ex>[r]^{{#2}} & \ {#4}
\ar@<0.4ex>[l]^{{#3}}
}
}

\def \mfs#1#2#3#4{\!
\xymatrix@C=1,5em{{#1} \! \ar@<0.2ex>[r]^{{#2}} & \! {#4}
\ar@<0.2ex>[l]^{{#3}}
}
\!}

\def \mfl#1#2#3#4{
\xymatrix@C=2.6em{{#1}\  \ar@<0.4ex>[r]^{{#2}} &\  {#4}
\ar@<0.2ex>[l]^{{#3}}
}
}

\def \mfss#1#2#3#4{\!
\xymatrix@C=1.5em{{#1} \ar@<0.3ex>[r]^{{#2}} & {#4}
\ar@<0.3ex>[l]^{{#3}}
}
\!}
\begin{document}
\title[Maximally graded m.f. for chain type polynomial]{Maximally-graded matrix factorizations for an invertible polynomial of chain type}
\date{\today}
\author{Daisuke Aramaki}
\address{Department of Mathematics, Graduate School of Science, Osaka University, 
Toyonaka Osaka, 560-0043, Japan}
\email{u935359a@ecs.osaka-u.ac.jp}
\author{Atsushi Takahashi}
\address{Department of Mathematics, Graduate School of Science, Osaka University, 
Toyonaka Osaka, 560-0043, Japan}
\email{takahashi@math.sci.osaka-u.ac.jp}
\begin{abstract}
In 1977, Orlik--Randell construct a nice integral basis of the middle homology group of 
the Milnor fiber associated to an invertible polynomial of chain type and they conjectured that it is represented by a distinguished basis of vanishing cycles.

The purpose of this paper is to prove the algebraic counterpart of the Orlik--Randell conjecture. 
Under the homological mirror symmetry, we may expect that the triangulated category of maximally-graded matrix factorizations 
for the Berglund--H\"{u}bsch transposed polynomial admits a full exceptional collection with a nice numerical property. 
Indeed, we show that the category admits a Lefschetz decomposition with respect to a polarization in the sense of Kuznetsov--Smirnov, 
whose Euler matrix are calculated in terms of the ``zeta function" of the inverse of the polarization. 

As a corollary, it turns out that the homological mirror symmetry holds at the level of lattices, 
namely, the Grothendieck group of the category with the Euler form is isomorphic to the middle homology group with the intersection form (with a suitable sign).
\end{abstract}
\maketitle
\section{Introduction}
Let $f$ be an invertible polynomial of chain type, namely, a polynomial of the form 
\begin{equation}
f= f(x_1,\dots, x_n):= x_1^{a_1}x_2 + x_2^{a_2}x_3 + \cdots + x_{n-1}^{a_{n-1}}x_n + x_n^{a_n}, \quad a_i\ge 2.
\end{equation}
One can associate an interesting algebro-geometric invariant, the triangulated category ${\rm HMF}^{L_f}_S(f)$ of $L_f$-graded matrix factorizations of $f$, 
where $S:=\CC[x_1,\dots, x_n]$ and $L_f$ is the maximal grading of $f$ (see Section 1).
This category ${\rm HMF}^{L_f}_S(f)$ is considered as an analogue of the bounded derived category of coherent sheaves on a smooth proper algebraic variety.
Let $\widetilde f$ be the Berglund--H\"{u}bsch transpose of $f$ defined by
\begin{equation}
\widetilde f=\widetilde f (x_1,\dots ,x_n):=x_1^{a_1} +x_1 x_2^{a_2}+\dots+x_{n-1}x_n^{a_n},
\end{equation}
which (also) has an isolated singularity only at the origin $0\in\CC^n$.
A distinguished basis of vanishing $(n-1)$-spheres 
in the Milnor fiber of $\widetilde f$ can be categorified to an $A_\infty$-category ${\rm Fuk}^{\to}(\widetilde f)$
called the directed Fukaya category. Its derived category $D^b{\rm Fuk}^\to(\widetilde f)$ is an invariant of the polynomial $\widetilde f$ as a triangulated category 
and $D^b{\rm Fuk}^\to(\widetilde f)$ has a full exceptional collection by construction.
The Berglund--H\"{u}bsch transposition of invertible polynomials together with the Berglund--Henningson duality of their symmetries \cite{BHe} yields several mirror symmetry conjectures. 
For example, the topological mirror symmetry conjecture is now proven in \cite{EGT}.
Here in this paper, we are interested in the homological mirror symmetry conjecture, namely, we expect the following equivalence of triangulated categories to hold:
\begin{conj}[{cf. \cite{ET,T2}}]
\begin{equation}
{\rm HMF}^{L_f}_S(f)\cong  D^b{\rm Fuk}^\to(\widetilde f).
\end{equation}
\end{conj}
There is much evidence of the above conjecture which follow from related results by several authors (cf. \cite{T1,KST1,KST2,Ue,FU1,FU2,LP}). 
\begin{rem}
In just a few days before the submission to the arXiv, a closely related work by Habermann--Smith \cite{HS} has appeared.
They seem to prove the conjecture for invertible polynomials in two variables, where in order to show the equivalence 
they seem to use the same full strongly exceptional collection as the one in an unpublished work by the second named author and Yoko Hirano (cf. \cite{T3}). 
\end{rem}
Based on this conjecture, it is natural to expect the existence of a full exceptional collection in ${\rm HMF}^{L_f}_S(f)$. 
Recently, Hirano--Ouchi \cite{HO} prove this by their semi-orthogonal decomposition theorem and an induction on the number of variables in $f$.
However, unfortunately, it was not clear how to choose matrix factorizations which form a full exceptional collection.
In 1977, Orlik--Randell \cite[Theorem~2.11]{OR} construct a nice $\ZZ$-basis of the middle homology group of the Milnor fiber associated to $\widetilde f$
and conjectured that the $\ZZ$-basis is represented by a distinguished basis of vanishing cycles \cite[Conjecture~4.1]{OR}.
The purpose of this paper is to prove the algebraic counterpart of the Orlik--Randell conjecture. 
More precisely, we explicitly construct a full exceptional collection motivated by their $\ZZ$-basis, which is moreover a Lefschetz decomposition in the sense of Kuznetsov--Smirnov \cite[Definition~2.3]{KuS}, 
whose Euler matrix are calculated in terms of the ``zeta function" of the inverse of the polarization. Most important part of our main theorem (Theorem~\ref{thm: main})
is the following 
\begin{thm}\label{thm: intro}
There exist an exceptional object $E\in {\rm HMF}^{L_f}_S(f)$ and an autoequivalence $\tau$ on ${\rm HMF}^{L_f}_S(f)$ such that 
$(E,\tau E,\dots ,\tau^{\widetilde \mu_n -1} E)$ is a full exceptional collection:
\begin{equation}
{\rm HMF}^{L_f}_S(f)\cong \langle E,\tau E,\dots ,\tau^{\widetilde \mu_n -1} E\rangle,
\end{equation}
whose Euler matrix is given by
\begin{equation}
\prod_{i=0}^{n}\left(1-N^{d_i}\right)^{(-1)^{n-i+1}},
\end{equation}
where $\widetilde \mu_n:=d_n-d_{n-1}+\dots +(-1)^{n-1}d_1+(-1)^{n}$, $d_i:=a_1\dots a_i$, is the Milnor number of $\widetilde f$,
and $N:=(\delta_{i+1,j})$ is the regular nilpotent $\widetilde \mu_n\times \widetilde \mu_n$-matrix.
\end{thm}
We give here an outline of this paper. In Section~2, after introducing the definition of invertible polynomials and their maximal gradings, we recall some definitions and facts
on matrix factorizations necessary in later sections.
Section~3 is the main part of this paper. 
We first show the existence of the polarization $\tau$ of index $\widetilde{\mu}$ in the sense of Kuzunetsov-Smirnov \cite[Definition~2.3]{KuS}, 
the autoequivalence $\tau$ on ${\rm HMF}^{L_f}_S(f)$ satisfying 
\begin{equation}
T^{-n-2k}\S = \tau^{-\widetilde{\mu}}, 
\end{equation}
for some $k\in\ZZ$ where $\S$ is the Serre functor on ${\rm HMF}^{L_f}_S(f)$.
This motivates us a Lefschetz decomposition with respect to $\tau$ due to \cite{HO}.
Note that the automorphism induced by $T^{-n}\S$ on the Grothendieck group can be considered as the mirror dual object to the Milnor monodromy 
associated to the singularity associated to $\widetilde f$. 
Therefore, we then show the same linear algebraic property as the Milnor monodromy have. 
Based on this observation, we find a candidate of the object $E$ in Theorem~\ref{thm: intro} and come to our main theorem in this paper (Theorem~\ref{thm: main}). 
Section~4 is devoted to the proof of Theorem~\ref{thm: main}.
It is Lemma~\ref{lem:key} that enables us to show that the given exceptional collection is full, whose $\ZZ$-graded version is also used in \cite{KST2}. 
It is also the key that one can calculate systematically the spaces of morphisms due to Dyckerhof (the $L_f$-graded modification of Proposition~4.3 in \cite{D}).
It is known by by Kreuzer--Skarke \cite{KrS} that any invertible polynomial is a Thom--Sebastiani sum of chain type polynomials and loop type polynomials.
For an invertible polynomial of loop type, two important problems are still open; 
to prove the existence of a full exceptional collection and to construct systematically a ``good one" among full exceptional collections.
We hope to solve these in the near future.
\bigskip
\noindent{\bf Acknowledgements.}
We would like to thank Yuuki Hirano and Genki Ouchi for valuable discussion.
The second named author is grateful to Wolfgang Ebeling who prompted him to look at the Orlik--Randell conjecture.
The second named author is supported by JSPS KAKENHI Grant Number JP16H06337.
\section{Preliminaries}
In this section, we set up several definitions which are used in the present paper. 
\subsection{Invertible polynomials}
Let $f(x_1,\ldots, x_n)\in \CC[x_1,\dots, x_n]$ be a non-degenerate weighted homogeneous polynomial, 
namely, a polynomial with an isolated singularity at the origin with the property that there are positive integers 
$w_1,\ldots ,w_n$ and $d$ such that 
$f(\lambda^{w_1} x_1, \ldots, \lambda^{w_n} x_n) = \lambda^d f(x_1,\ldots ,x_n)$ 
for $\lambda \in \CC^\ast$. 
\begin{defn}
A non-degenerate weighted homogeneous polynomial $f(x_1,\ldots ,x_n)$ is called {\em invertible} if 
the following conditions are satisfied:
\begin{enumerate}
\item the number of variables ($=n$) coincides with the number of monomials 
in the polynomial $f(x_1,\ldots x_n)$, 
namely, 
\[
f(x_1,\ldots ,x_n)=\sum_{i=1}^n c_i\prod_{j=1}^nx_j^{E_{ij}}
\]
for some coefficients $c_i\in\CC^\ast$ and non-negative integers 
$E_{ij}$ for $i,j=1,\ldots, n$,
\item the matrix $E:=(E_{ij})$ is invertible over $\QQ$.
\end{enumerate}
\end{defn}
Since we work over the complex number field $\CC$, without loss of generality one may assume that $c_i=1$ for $i=1,\ldots, n$, 
which can be achieved by rescaling the variables.
There is also a classification of invertible polynomials. It is known by \cite{KrS} that an invertible polynomial $f$ is a Thom--Sebastiani sum of invertible polynomials of the following types:
\begin{enumerate}
\item[1)] $x_1^{a_1}x_2 + x_2^{a_2}x_3 + \cdots + x_{n-1}^{a_{n-1}}x_n + x_n^{a_n}$ (chain type; $n\ge 1$);
\item[2)] $x_1^{a_1}x_2 + x_2^{a_2}x_3 + \cdots + x_{n-1}^{a_{n-1}}x_n + x_n^{a_n}x_1$ (loop type; $n\ge 2$).
\end{enumerate}
\begin{defn}
Let $f(x_1,\dots ,x_n)=\displaystyle\sum_{i=1}^na_i\prod_{j=1}^nx_j^{E_{ij}}$ be an invertible polynomial.
Consider the free abelian group $\displaystyle\bigoplus_{i=1}^n\ZZ\vec{x}_i\oplus \ZZ\vec{f}$ 
generated by the symbols $\vec{x}_i$ for the variables $x_i$ for $i=1,\dots, n$
and the symbol $\vec{f}$ for the polynomial $f$.
The {\em maximal grading} $L_f$ of the invertible polynomial $f$ 
is the abelian group defined by the quotient 
\begin{equation}
L_f:=\left.\left(\bigoplus_{i=1}^n\ZZ\vec{x}_i\oplus \ZZ\vec{f}\right)\right/\left(\vec{f}-\sum_{j=1}^nE_{ij}\vec{x}_j;\ i=1,\dots ,n\right).
\end{equation}
Note that $L_f$ is an abelian group of rank $1$ which is not necessarily free.
\end{defn}
\subsection{Matrix factorizations}
In this subsection, we recall some properties of $L_f$-graded matrix factorizations attached to an invertible polynomial $f$. 
All the results are slight modifications of those in \cite{KST1, KST2, Or1, Or2}.
Set $S:=\CC[x_1,\dots. x_n]$, which is naturally an $L_f$-graded $\CC$-algebra;
\begin{equation}
S=\bigoplus_{\vec{l}\in L_f} S_{\vec{l}}.
\end{equation}
For any $L_f$-graded $S$-module $M$, $M(\vec{l})$ denotes 
the grading shift by $\vec{l}\in L_f$ of $M;$
\begin{equation}
M(\vec{l})=\bigoplus_{\vec{l'}\in L_f} M(\vec{l})_{\vec{l'}},\quad M(\vec{l})_{\vec{l'}}:=M_{\vec{l}+\vec{l'}},
\end{equation}
which induces the autoequivalence functor $(\vec{l})$ for $\vec{l}\in L_f$ on ${\rm gr}^{L_f}\text{-}S$, 
the category of finitely generated $L_f$-graded $S$-modules. 
\begin{defn}[Eisenbud \cite{Ei}]
Let $F_0,F_1$ be $L_f$-graded free modules and $f_0:F_0\to F_1$, $f_1:F_1\to F_0(\vec{f})$ be 
morphisms such that $f_1\circ f_0=f\cdot{\rm id}_{F_0}, f_0(\vec{f})\circ f_1=f\cdot{\rm id}_{F_1}$.
The tuple $(F_0,F_1,f_0,f_1)$ is called a $L_f$-graded {\it matrix factorization} of $f$
and denoted by 
\begin{equation}
\overline{F}:=\Big(\mf{F_0}{f_0}{f_1}{F_1}\Big).
\end{equation}
\end{defn}
For an $L_f$-graded matrix factorization $\overline{F}$, 
the rank of $F_0$ must coincide with the one of $F_1$, 
which we shall call the {\em size} of the matrix factorization $\overline{F}$.
\begin{defn}
Let $\overline{F}=\Big(\mf{F_0}{f_0}{f_1}{F_1}\Big)$ and 
$\overline{F'}=\Big(\mf{F'_0}{f'_0}{f'_1}{F'_1}\Big)$ be $L_f$-graded matrix factorizations.
\begin{enumerate}
\item 
A {\it morphism} $\phi: \overline{F}\to \overline{F'}$ of $L_f$-graded matrix factorizations 
from $\overline{F}$ to $\overline{F'}$
is a pair $\phi=(\phi_0,\phi_1)$ where $\phi_0:F_0\to F'_0$ and $\phi_1:F_1\to F'_1$ are morphisms in ${\rm gr}^{L_f}\text{-}S$ such that 
$\phi_1\circ f_0=f'_0\circ \phi_0$ and $\phi_0(\vec{f})\circ f_1=f'_1\circ \phi_1$.
\item 
The morphism $\phi=(\phi_0,\phi_1):\overline{F}\to \overline{F'}$ 
is called {\it null-homotopic} if there is a pair $\psi=(\psi_0,\psi_1)$ where $\psi_0:F_0\to F'_1(-\vec{f})$ and $\psi_1:F_1\to F'_0$ are 
morphisms in ${\rm gr}^{L_f}\text{-}S$ such that 
$\phi_0=f'_1(-\vec{f})\circ \psi_0+\psi_1\circ f_0$ and 
$\phi_1=\psi_0(\vec{f})\circ f_1+f'_0\circ \psi_1$.
\end{enumerate}
\end{defn}
\begin{prop}
The category ${\rm MF}^{L_f}_S(f)$ of $L_f$-graded matrix factorizations of $f$ 
is a Frobenius category whose morphisms factoring through projectives coincide with 
null-homotopic morphisms.
Therefore, its stable category 
\begin{equation}
{\rm HMF}^{L_f}_S(f):=\underline{{\rm MF}}^{L_f}_S(f)
\end{equation}
has a natural structure of a triangulated category due to \cite{Ha}.
\qed
\end{prop}
We shall denote by $T$ the translation functor on ${\rm HMF}^{L_f}_S(f)$. 
By definition of the triangulated structure on ${\rm HMF}^{L_f}_S(f)$, 
one easily see the following. 
\begin{prop}
Each exact triangle in ${\rm HMF}^{L_f}_S(f)$ is isomorphic to a triangle of the form 
\begin{equation}
\overline{F}\stackrel{[\phi]}{\longrightarrow}
\overline{F'}\stackrel{}{\longrightarrow}
C(\phi)\stackrel{}{\longrightarrow}
T\overline{F}
\end{equation}
for some $\overline{F},\overline{F'}\in {\rm HMF}^{L_f}_S(f)$ 
and $\phi\in{\rm MF}^{L_f}_S(f)(\overline{F},\overline{F'})$, a lift of $[\phi]\in{\rm HMF}^{L_f}_S(f)(\overline{F},\overline{F'})$. 
\qed
\end{prop}
Since $f$ has an isolated singularity at the origin, we have the following
\begin{prop}
The category ${\rm HMF}^{L_f}_S(f)$ has the following properties:
\begin{enumerate}
\item It is finite. Namely, for all $\overline{F}, \overline{F'}\in {\rm HMF}^{L_f}_S(f)$, we have
\begin{equation}
\sum_{p\in\ZZ}\dim_\CC {\rm HMF}^{L_f}_S(f)(\overline{F},T^p \overline{F'})<\infty.
\end{equation}
\item 
It is idempotent complete. Namely, 
for any $\overline{F}\in {\rm HMF}^{L_f}_S(f)$ and any idempotent $e\in {\rm HMF}^{L_f}_S(f)(\overline{F},\overline{F})$, 
there exists an object $\overline{F'}\in {\rm HMF}^{L_f}_S(f)$ and a pair of morphisms $\phi\in {\rm HMF}^{L_f}_S(f)(\overline{F},\overline{F'})$ and $\phi'\in {\rm HMF}^{L_f}_S(f)(\overline{F'},\overline{F})$ such that 
$\phi'\circ \phi=e$ and $\phi\circ \phi'={\rm id}_{\overline{F'}}$
\end{enumerate}
\qed
\end{prop}
The autoequivalence functor $(\vec{l})$ for $\vec{l}\in L_f$ on ${\rm gr}^{L_f}\text{-}S$ 
induces an autoequivalence on ${\rm HMF}^{L_f}_S(f)$, 
which we denote by the same symbol $(\vec{l})$.
Explicitly, the action of $(\vec{l})$ takes an object $\overline{F}$ to the object
\begin{equation}
\overline{F}(\vec{l}):=
\Big(\mf{F_0(\vec{l})}{f_0(\vec{l})}{f_1(\vec{l})}{F_1(\vec{l})}\Big),
\end{equation}
and takes a morphism $\phi=(\phi_0,\phi_1)$ to the morphism 
$\phi(\vec{l}):=(\phi_0(\vec{l}),\phi_1(\vec{l}))$.
Similary, the translation functor $T$ on ${\rm HMF}^{L_f}_S(f)$ takes an object 
$\overline{F}$ to the object
\begin{equation}
T\overline{F}
:=\Big(\mf{F_1}{-f_1\quad}{-f_0(\vec{f})\quad}{F_0(\vec{f})}\Big),
\end{equation}
and takes a morphism $\phi=(\phi_0,\phi_1)$ to the morphism 
$T(\phi):=(\phi_1,\phi_0(\vec{f}))$.
From these descriptions of autoequivalences $(\vec{l})$ and $T$, 
we obtain the following:
\begin{prop}
On the category ${\rm HMF}^{L_f}_S(f)$, we have $T^2=(\vec{f})$.
\qed
\end{prop}
It is also known that Auslander--Reiten duality \cite{AR} yields the existence of the Serre functor on ${\rm HMF}^{L_f}_S(f)$.
\begin{prop}
The functor 
\begin{equation}
\S:=T^{n-2}\circ (-\vec{\varepsilon}_f)=T^n\circ(-\vec{x}_1-\dots -\vec{x}_n),\quad \vec{\varepsilon}_f:=\left(\sum_{i=1}^n\vec{x}_i\right) -\vec{f},
\end{equation}
defines the Serre functor on  ${\rm HMF}^{L_f}_S(f)$, namely, there exists a bi-functorial isomorphism
\begin{equation}
{\rm HMF}^{L_f}_S(f)(\overline{F}, \overline{F'})\cong 
 {\rm HMF}^{L_f}_S(f)(\overline{F'},\S\overline{F})^\ast,
\end{equation}
where ${}^\ast$ denotes the duality over $\CC$.
\qed
\end{prop}
We also recall some notions necessary for statements below. 
\begin{defn}
Let $\T$ be a $\CC$-linear triangulated category $\T$ with the translation functor $T$.
\begin{enumerate}
\item
An object $E$ in $\T$ is called an {\em exceptional object} (or is called {\em exceptional}) if 
$\T(E,E)=\CC\cdot {\rm id}_E$ and $\T(E,T^p E)=0$ when $p\ne 0$.
\item 
An {\em exceptional collection} $\E=(E_1,\dots, E_\mu)$ in $\T$ is a finite set of exceptional objects  
satisfying the condition $\T(E_i,T^p E_j)=0$ for all $p$ and $i>j$.
\item 
An exceptional collection $\E=(E_1,\dots, E_\mu)$ in $\T$ is called 
a {\em strongly exceptional collection} if $\T(E_i,T^p E_j)=0$ for all $p\ne 0$ and $i,j=1,\dots, \mu$.
\item 
An exceptional collection $\E=(E_1,\dots, E_\mu)$ in $\T$ is called {\em full}
if the smallest full triangulated subcategory of $\T$ containing all elements in $\E$
is equivalent to $\T$ as a triangulated category.
\end{enumerate}
\end{defn}
\begin{defn}[cf. {\cite[Section~2.3]{D}}]
Let $M$ be an $L_f$-graded $S$-module of the form $S/(p_1, \dots, p_s)$ such that 
$p_i\in S_{\vec{p}_i}$ for some $\vec{p}_i\in L_f$, $p_1,\dots , p_s$ forms a regular sequence and $f\in (p_1, \dots, p_s)$.
Write $f=p_1h_1+\dots +p_s h_s$ with $h_i\in S_{\vec{f}-\vec{p}_i}$ and put $P:=\oplus_{i=1}^{s}S(-\vec{p}_i)$.
Then the Koszul resolution of $M$ as an $L_f$-graded $S$-module 
\begin{equation}
0\longrightarrow \wedge^s P\longrightarrow \wedge^{s-1} P \longrightarrow \dots \longrightarrow \wedge^1 P\longrightarrow \wedge^{0} P=S\longrightarrow
M\longrightarrow 0,
\end{equation}
yields the $L_f$-graded matrix factorization $\overline{F}=(\mfs{F_0}{f_0}{f_1}{F_1})$ of $f$ such that
\begin{equation}
F_0:=\bigoplus_{k}\left(\wedge^{2k} P\right)(k\vec{f}),\quad F_1:=\bigoplus_{k}\left(\wedge^{2k-1} P\right)(k\vec{f}),
\end{equation}
which is unique up to isomorphism in ${\rm HMF}^{L_f}_{S}(f)$ and is called the {\em stabilization} of $M$ and will be denoted by $M^{stab}$.
\end{defn}
For matrix factorizations from stabilizations, there is a convenient way to calculate their spaces of morphisms. See \cite[Proposition~4.3]{D}.  
Note that the category ${\rm HMF}^{L_f}_{S}(f)$ always contains an object $\CC^{stab}$, the stabilization of the simple $L_f$-graded $S$-module 
$\CC=S/(x_1,\dots, x_n)$, which plays an important role. 
In particular, we have the following key lemma, the $L_f$-graded version of the one by Kajiura-Saito-Takahashi \cite[Theorem~4.5]{KST2}.
\begin{lem}\label{lem:key}
Let $\E$ be an exceptional collection in ${\rm HMF}^{L_f}_{S}(f)$ and 
$\langle \E\rangle$ be the smallest full triangulated subcategory of ${\rm HMF}^{L_f}_{S}(f)$ containing all elements in $\E$.
Suppose that $\langle \E\rangle$ satisfies the following conditions$:$
\begin{enumerate}
\item $\langle \E\rangle$ is closed under the grading shift $(\vec{l})$ for all $\vec{l}\in L_f$.
\item $\langle \E\rangle$ contains $\CC^{stab}$.
\end{enumerate}
Then, the natural fully faithful functor $\langle \E\rangle\longrightarrow {\rm HMF}^{L_f}_{S}(f)$ is a triangulated equivalence.
\qed
\end{lem}
It is natural from the homological mirror symmetry conjecture to expect that the category ${\rm HMF}^{L_f}_S(f)$ admits 
a full exceptional collection. Indeed, there are many results which support this expectation (cf. \cite{T1, KST1,KST2, Ue, FU1, FU2, LP}).
For $f$ of chain type, recently Hirano--Ouchi \cite{HO} prove the existence of the full exceptional collection.
The main result in this paper is that, furthermore, there is a ``good choice'' among full exceptional collections from both combinatorial and geometric points of view. 
\section{Main result}
From now on, we shall only consider an invertible polynomial $f$ of chain type;
\[
f=x_1^{a_1}x_2 + x_2^{a_2}x_3 + \cdots + x_{n-1}^{a_{n-1}}x_n + x_n^{a_n}.
\]
For simplicity, we assume that $a_i\ge 2$ for all $i=1,\dots, n$.
Set 
\begin{equation}
d_i:=a_1\cdots a_i,\quad d_0:=1,
\end{equation}
and define $\widetilde \mu_i$ inductively by
\begin{equation}
\widetilde \mu_{i}:=d_{i}-\widetilde \mu_{i-1},\quad \widetilde \mu_0:=1.
\end{equation}
It is important that $L_f/\ZZ\vec{f}$ is generated by $\vec{x}_1$, which follows from the relations in $L_f$;
\begin{equation}
\vec{f}=a_1\vec{x}_1+\vec{x_2}=a_2\vec{x}_2+\vec{x_3}=\dots =a_{n-1}\vec{x}_{n-1}+\vec{x_n}=a_n\vec{x}_n.
\end{equation}
Indeed, $\vec{x}_i=(-1)^{i-1}d_{i-1}\vec{x}_1$ for $i=2,\dots, n$ and $d_n\vec{x}_1=\vec{0}$ in $L_f/\ZZ\vec{f}$.
In particular, from the description of the Serre functor $\S$ on ${\rm HMF}^{L_f}_{S}(f)$, this implies the following  
\begin{prop}
There exists an integer $k\in\ZZ$ such that 
\begin{equation}
T^{-n-2k}\S=
\begin{cases}
(\widetilde \mu_n \vec{x}_1) & \text{if}\quad n=2m+1,\ m\in \ZZ_{\ge 0},\\
(-\widetilde \mu_n \vec{x}_1) & \text{if}\quad n=2m,\ m\in\ZZ_{\ge 1}.
\end{cases}
\end{equation}
\qed
\end{prop}
It is shown in \cite[Corollary~4.6]{HO} that ${\rm HMF}^{L_f}_{S}(f)$ has a full exceptional collection $\E=(E_1,\dots,E_{\widetilde \mu_n})$.
The above proposition means that, if $n$ is odd (resp. even), $\tau:=(-\vec{x}_1)$ (resp. $\tau:=(\vec{x}_1)$) is a polarization of ${\rm HMF}^{L_f}_{S}(f)$ of index $\widetilde \mu_n$
in the sense of Kuznetsov--Smirnov \cite[Definition~2.1]{KuS}.
Therefore, it is natural to expect the existence of an exceptional object $E$ such that 
\[
{\rm HMF}^{L_f}_{S}(f)\cong \langle E, \tau E, \dots, \tau^{\widetilde\mu_n-1} E\rangle,
\]
a Lefschetz decomposition of ${\rm HMF}^{L_f}_{S}(f)$ with respect to the above polarization $\tau$
in the sense of \cite[Definition~2.3]{KuS}.
In order to specify the numerical properties to be satisfied by this exceptional object $E$, 
following \cite{OR} we first introduce a polynomial $\varphi_n(t)$ in $t$ of degree $\widetilde \mu_n$;
\begin{equation}
\varphi_n(t):=\prod_{i=0}^{n}\left(1-t^{d_i}\right)^{(-1)^{n-i}}.
\end{equation}
It will turn out later that $\varphi_n(t)$ is the ``zeta function" of the automorphism induced by $\tau^{-1}$ 
on the Grothendieck group $K_0({\rm HMF}^{L_f}_{S}(f))$.
\begin{prop}\label{prop:M_1}
Let $c'_i$ be the coefficients of $t^i$ in $\varphi_{n}(t)$. 
\begin{enumerate}
\item
We have 
\begin{equation}
c'_0=1,\quad c'_{\widetilde \mu_n-i}=(-1)^{n+1}c'_i.
\end{equation}
\item
We have 
\begin{equation}
\varphi_n(t)=\det ({\bf 1}-tM_1),
\end{equation}
where $M_1$ is the following $\widetilde \mu_n\times \widetilde \mu_n$-matrix
\begin{equation}
M_1:=
\begin{pmatrix}
-c'_1 & 1 & 0 & \cdots & 0 \\
-c'_2 & 0 & 1 & \ddots & \vdots\\
\vdots & 0 & \ddots & \ddots & 0\\
\vdots & \vdots & \ddots & \ddots & 1 \\
-c'_{\widetilde \mu_n} & 0 & \cdots &0 & 0
\end{pmatrix}.
\end{equation}
\end{enumerate}\qed
\end{prop}
Define a $\widetilde \mu_n\times \widetilde \mu_n$-matrix $\chi_n$ by 
\begin{equation}\label{eq:inverse}
\chi_n:=\frac{1}{\varphi_n(N)}=\prod_{i=0}^{n}\left(1-N^{d_i}\right)^{(-1)^{n-i+1}},
\end{equation}
where $N$ is the following nilpotent $\widetilde \mu_n\times \widetilde \mu_n$-matrix
\begin{equation}
N:=
\begin{pmatrix}
0 & 1 & 0 & \cdots & 0 \\
\vdots & 0 & 1 & \ddots & \vdots\\
\vdots & \ddots & \ddots & \ddots & 0\\
\vdots & \ddots & \ddots & \ddots & 1 \\
0 & \cdots & \cdots &\cdots & 0
\end{pmatrix}.
\end{equation}
More concretely, $\chi_n$ is given by 
\begin{equation}
\chi_n=
\begin{pmatrix}
c_0 & c_1 & c_2 & \cdots & \cdots & c_{\widetilde \mu_n -2} & c_{\widetilde \mu_n-1}\\
0 & c_0 & c_1 & c_2 & \ddots & \ddots & c_{\widetilde \mu_n -2}\\
0 & 0 & c_0 & c_1 & c_2 & \ddots & \vdots\\
\vdots & \ddots & \ddots & \ddots & \ddots & \ddots & \vdots \\
\vdots & \ddots & \ddots &\ddots &c_0 & c_1& c_2\\
\vdots & \ddots & \ddots &\ddots &0 & c_0& c_1\\
0 & \cdots & \cdots & \cdots & 0 & 0 &c_0
\end{pmatrix},
\end{equation}
where $c_i$ is the coefficient of $t^i$ in $1/\varphi_n(t)$; 
\begin{equation}
\frac{1}{\varphi_n(t)}=c_0+c_1 t+\dots +c_{\widetilde \mu_n-1} t^{\widetilde \mu_n-1}+c_{\widetilde \mu_n}t^{\widetilde \mu_n}+\cdots.
\end{equation}
Note that $c_0=1$ and for each positive integer $j$ we have
\begin{equation}\label{eq: inverse matrix}
\sum_{i=0}^\infty c_i c'_{j-i}=0,
\end{equation}
where we put $c'_{k}=0$ if $k<0$.
Set 
\begin{equation}
M:=(-1)^n\chi_n^{-1}\chi_n^T,
\end{equation}
and consider its ``zeta function"
\begin{equation}
\Phi_n(t):={\rm det}\left({\bf 1}-tM\right).
\end{equation}
Note that $M$ will be considered as a matrix representation of the automorphism induced by $T^{-n}\S$ 
on the Grothendieck group $K_0({\rm HMF}^{L_f}_{S}(f))$.
\begin{prop}\label{prop:M}
Let the notations be as above.
\begin{enumerate}
\item
We have
\begin{equation}
M=M_1^{\widetilde \mu_n}.
\end{equation}
\item
We have 
\begin{equation}
\Phi_n(t)=\prod_{i=0}^{n}\left(1-t^{\frac{d_i}{e_i}}\right)^{(-1)^{n-i}e_i},\quad e_i:={\rm gcd}(d_i,\widetilde \mu_n).
\end{equation}
In particular, $t^{\widetilde \mu_n}\cdot \Phi_n(t^{-1})$ is the characteristic polynomial of 
the Milnor monodromy of the singularity associated to 
the Berglund--H\"{u}bsch transpose $\widetilde f$ of $f$ \cite{BHu}:
\begin{equation}
\widetilde f=\widetilde f (x_1,\dots ,x_n)=x_1^{a_1} +x_1 x_2^{a_2}+\dots+x_{n-1}x_n^{a_n}.
\end{equation}
\end{enumerate}
\end{prop}
\begin{pf}
Using Proposition~\ref{prop:M_1} (i) and \eqref{eq: inverse matrix}, we obtain by a direct calculation that
\[
\chi_n M_1^{\widetilde \mu_n}=\left(\left(\dots\left(\left(\chi_n M_1\right)M_1\right)\dots \right)M_1\right)=(-1)^n\chi_n^T,
\]
which implies (i). 
Together with (i), Proposition~\ref{prop:M_1} (ii) yields the first statement of (ii).
The characteristic polynomial of the Milnor monodromy of the singularity associated to $\widetilde f$ 
can be computed using Varchenko's method \cite[Theorem~4.1.3]{Va}, 
which coincides with $t^{\widetilde \mu_n}\cdot \Phi_n(t^{-1})$.
\qed
\end{pf}
\begin{rem}
In order to have Proposition~\ref{prop:M} (i), one only needs a polynomial satisfying Proposition~\ref{prop:M_1} (i) for some $n$
and $1/\varphi_n(t)$ may not be necessarily a polynomial.
For example, one can apply the same proof for $(1-t)^{n+1}$, the ``zeta function" of the automorphism induced by $\O(-1)$ 
on the Grothendieck group of the bounded derived category $\D^b(\PP^n)$ 
of coherent sheaves on the projective space $\PP^n$, 
which admits a Lefschetz decomposition $\D^b(\PP^n)\cong \langle \O(0), \O(1),\dots, \O(n)\rangle$ with the polarization $\O(1)$.
Note that $(1-t)^{-n-1}$ is the Poincar\'e series of the ($\ZZ$-graded) polynomial ring in $n+1$-variables, which also calculates 
the dimension of ${\rm Hom}(\O(0),\O(i))$, $i=0,\dots, n$.
\end{rem}
\begin{rem}
Equivalent results are obtained in Balnojan-Hertling \cite{BaHe}, in which paper they consider these in order to reconstruct 
spectral numbers of a singularity starting from an upper triangular matrix, its Seifert form.
According to \cite{BaHe}, Horocholyn \cite{Ho} seems to be the first one who considers such $\widetilde \mu_n$-th root as $M_1$ here. 
As is already comment in \cite{BaHe}, Orlik--Randell should have known these including its geometric meaning although they did not write properly in \cite{OR}.
As far as we know, the use of the ``zeta functions" and the equation \eqref{eq:inverse}, which simplified the proof so much, 
and a natural algebraic explanation behind these facts are new.
\end{rem}
Now we can state our main theorem in this paper, which is motivated by the homological mirror symmetry and 
the conjecture by Orlik--Randell \cite[Conjecture~4.1]{OR} in singularity theory. In 1977, Orlik--Randell 
found a $\ZZ$-basis of the middle homology group $H_{n-1}(\widetilde f^{-1}(1),\ZZ)$ of the Milnor fiber $\widetilde f^{-1}(1)$ 
in which the intersection form and the Milnor monodromy are represented by the matrices
$(-1)^{(n-1)(n-2)/2}(\chi_n^{-1} +(\chi_n^{-1})^T)$ and $M^{-1}$. Then they
conjectured that the $\ZZ$-basis is represented by a distinguished basis of vanishing cycles, whose Seifert matrix is given by $(-1)^{(n-1)(n-2)/2}\chi_n^{-1}$.
\begin{thm}\label{thm: main}
Let the notations be as above.
Set 
\begin{equation}
E_i:=\begin{cases}
\left(S/(x_1, x_3, \dots, x_{2m+1})\right)^{stab}(-i\vec{x}_1) & \text{if}\quad n=2m+1,\ m\in \ZZ_{\ge 0},\\
\left(S/(x_2,x_4, \dots, x_{2m})\right)^{stab}(i\vec{x}_1) & \text{if}\quad n=2m,\ m\in\ZZ_{\ge 1}.
\end{cases}
\end{equation}
Then, for each $a\in \ZZ$, $(E_a,\dots E_{\widetilde \mu_n-1+a})$ forms a full exceptional collection in ${\rm HMF}^{L_f}_S(f)$.
Moreover, its Euler matrix coincides with $\chi_n$:
\begin{equation}
\chi_n=(\chi(E_i,E_j))_{i,j=a}^{\widetilde \mu_n-1+a},\quad \chi(E_i,E_j):=\sum_{p\in\ZZ}(-1)^p\dim_\CC {\rm HMF}^{L_f}_S(f)(E_i,T^pE_j).
\end{equation}
\end{thm}
\begin{rem}
Once we have a full exceptional collection whose Euler matrix is given by $\chi_n$, then we also have another one with the Euler matrix 
$\chi_n^{-1}$. 
In this sense, Theorem~\ref{thm: main} solves the homological mirror symmetric dual statement of the Orlik--Randell's conjecture.
\end{rem}
We immediately have the following 
\begin{cor}
For an invertible polynomial of chain type, the homological mirror symmetry at the level of the Grothendieck group holds.
Namely, we have the following isomorphism of lattices 
\begin{equation}
\left(K_0({\rm HMF}^{L_f}_S(f)), \chi+\chi^T\right)\cong \left(H_{n-1}(\widetilde f^{-1}(1),\ZZ),(-1)^{\frac{(n-1)(n-2)}{2}}I\right),
\end{equation}
where $\chi$ is the Euler form and $I$ is the intersection form.
\qed
\end{cor}
For some special cases, Theorem~\ref{thm: main} is already known. For example, for $f=x_1^{a_1}$ it follows from \cite{T1}, and for $f=x_1^{2k}x_2+x_2^2$ ($k\ge 1$) and 
$f=x_1^2x_2+x_2^{a_2}x_3+x_3^2$ one can derive the above statement from the corresponding result in \cite{KST1}.
If $n=2$, then the full exceptional collection in Theorem~\ref{thm: main} is strong, whose endomorphism algebra 
is a Nakayama algebra.
\begin{cor}\label{thm:n=2}
For $n=2$, we have an equivalence of triangulated categories 
\begin{equation}
{\rm HMF}^{L_f}_S(f)\cong \D^b{\rm mod}\text{-} A_{\widetilde \mu_2}(a_1),
\end{equation}
where $A_{\widetilde \mu_2}(a_1)$ is the Nakayama algebra given by 
the equi-oriented $A_{\widetilde \mu_2}$-quiver 
\begin{equation}
1 \stackrel{x}{\longrightarrow} 2 \stackrel{x}{\longrightarrow} 3\stackrel{x}{\longrightarrow} \dots 
\stackrel{x}{\longrightarrow} \widetilde \mu_2-1 \stackrel{x}{\longrightarrow} \widetilde \mu_2,
\end{equation}
with all relations $x^{a_1}=0$. 
\end{cor}
\begin{pf}
This is a direct consequence of Theorem~\ref{thm: main} and Lemma~\ref{lem:4.1} below for $n=2$.
\qed
\end{pf}
\section{Proof}
We will apply our key lemma, Lemma~\ref{lem:key}, to prove our main theorem. 
\subsection{Morphisms}
Let the notations be as in the previous sections.
Since the object $E_i$ is the stabilization of an $L_f$-graded $S$-module, 
we obtain the following two lemmas by a direct calculation.
\begin{lem}\label{lem:4.1}
Suppose that $n=2m$, $m\in \ZZ_{\ge 1}$.
\begin{enumerate}
\item
For each $i\in \ZZ$, we have a natural isomorphism of $L_f$-graded $\CC$-algebras
\[
\bigoplus_{\vec{l}\in L_f}{\rm HMF}^{L_f}_{S}(f)(E_i, E_i(\vec{l}))
\cong S/(x_2, x_4,\dots, x_{2m}, x_1^{a_1},x_3^{a_3},\dots, x_{2m-1}^{a_{2m-1}}),
\]
where the product of the LHS is the one induced by the composition map
\begin{eqnarray*}
& &{\rm HMF}^{L_f}_{S}(f)(E_i, E_i(\vec{l}))\times {\rm HMF}^{L_f}_{S}(f)(E_i, E_i(\vec{l}'))\\
&\cong &{\rm HMF}^{L_f}_{S}(f)(E_i(\vec{l}'), E_i(\vec{l}+\vec{l}'))\times {\rm HMF}^{L_f}_{S}(f)(E_i, E_i(\vec{l}))\\
&\longrightarrow &
{\rm HMF}^{L_f}_{S}(f)(E_i, E_i(\vec{l}+\vec{l}')).
\end{eqnarray*}
\item
For each $i\in \ZZ$, we have
\[
\bigoplus_{\vec{l}\in L_f}{\rm HMF}^{L_f}_{S}(f)(E_i, TE_i(\vec{l}))=0.
\]
\end{enumerate}
\end{lem}
\begin{pf}
We can write $f$ as 
\[
f=x_2(x_1^{a_1}+x_2^{a_2-1}x_3)+x_4(x_3^{a_3}+x_4^{a_4-1}x_5)+\dots +x_{2m}(x_{2m-1}^{a_{2m-1}}+x_{2m}^{a_{2m}-1})
\]
to describe the stabilization of $S/(x_2, x_4,\dots, x_{2m})$.
By this expression, all the statements are straightforward since $a_i\ge 2$ for all $i=1,\dots, 2m$.
\qed
\end{pf}
\begin{lem}\label{lem:4.2}
Suppose that $n=2m+1$, $m\in \ZZ_{\ge 0}$.
\begin{enumerate}
\item
For each $i\in \ZZ$, we have a natural isomorphism of $L_f$-graded $\CC$-algebras
\[
\bigoplus_{\vec{l}\in L_f}{\rm HMF}^{L_f}_{S}(f)(E_i, E_i(\vec{l}))
\cong S/(x_1, x_3,\dots, x_{2m+1}, x_2^{a_2},x_4^{a_4},\dots, x_{2m}^{a_{2m}}),
\]
\item
For each $i\in \ZZ$, 
$\oplus_{\vec{l}\in L_f}{\rm HMF}^{L_f}_{S}(f)(E_i, TE_i(\vec{l}))$
has a structure of an $L_f$-graded module over $\oplus_{\vec{l}\in L_f}{\rm HMF}^{L_f}_{S}(f)(E_i, E_i(\vec{l}))$ 
by the composition map
\begin{eqnarray*}
& &{\rm HMF}^{L_f}_{S}(f)(E_i, TE_i(\vec{l}))\times {\rm HMF}^{L_f}_{S}(f)(E_i, E_i(\vec{l}'))\\
&\cong &{\rm HMF}^{L_f}_{S}(f)(E_i(\vec{l}'), TE_i(\vec{l}+\vec{l}'))\times {\rm HMF}^{L_f}_{S}(f)(E_i, E_i(\vec{l}))\\
&\longrightarrow &
{\rm HMF}^{L_f}_{S}(f)(E_i, TE_i(\vec{l}+\vec{l}')), 
\end{eqnarray*}
which is free of rank one generated by an element in ${\rm HMF}^{L_f}_{S}(f)(E_i, TE_i(-\vec{x_1}))$. 
Moreover, the composition map
\begin{eqnarray*}
& &{\rm HMF}^{L_f}_{S}(f)(E_i, TE_i(\vec{l}))\times {\rm HMF}^{L_f}_{S}(f)(E_i, TE_i(\vec{l}'))\\
&\cong &{\rm HMF}^{L_f}_{S}(f)(TE_i(\vec{l}'), T^2E_i(\vec{l}+\vec{l}'))\times {\rm HMF}^{L_f}_{S}(f)(E_i, TE_i(\vec{l}))\\
&\longrightarrow &
{\rm HMF}^{L_f}_{S}(f)(E_i, E_i(\vec{f}+\vec{l}+\vec{l}'))
\end{eqnarray*}
is zero for all $\vec{l},\vec{l}'\in L_f$.
\end{enumerate}
\end{lem}
\begin{pf}
We can write $f$ as 
\[
f=x_1(x_1^{a_1}x_2)+x_3(x_2^{a_2}+x_3^{a_3-1}x_4)+x_5(x_4^{a_4}+x_5^{a_5-1}x_6)+\dots +x_{2m+1}(x_{2m}^{a_{2m}}+x_{2m+1}^{a_{2m+1}-1})
\]
to describe the stabilization of $S/(x_1, x_3,\dots, x_{2m+1})$.
By this expression, all the statements are straightforward since $a_i\ge 2$ for all $i=1,\dots, 2m+1$.
\qed
\end{pf}
The following is a direct consequence of the above lemmas.
\begin{cor}
For each $a\in\ZZ$, $(E_a,\dots E_{\widetilde \mu_n-1+a})$ forms an exceptional collection.
\qed
\end{cor}
Since $\vec{x}_i= (-1)^{i-1}d_{i-1}\vec{x}_1$ in $L_f/\ZZ\vec{f}$ and $(\vec{f})=T^2$, by calculating the generating function of 
dimensions of spaces of morphisms, we obtain the following corollaries.
\begin{cor}
Suppose that $n=2m$, $m\in \ZZ_{\ge 1}$. 
We have 
\begin{eqnarray*}
(\chi(E_i,E_j))_{i,j=a}^{\widetilde \mu_n-1+a}&=&
\frac{1-N^{a_1}}{1-N}\cdot \frac{1-(N^{d_2})^{a_3}}{1-N^{d_2}}\cdot \dots \cdot \frac{1-(N^{d_{2m-2}})^{a_{2m-1}}}{1-N^{d_{2m-2}}}\\
&=&\prod_{i=0}^{n}\left(1-N^{d_i}\right)^{(-1)^{n-i+1}}=\frac{1}{\varphi_n(N)}=\chi_n.
\end{eqnarray*}
\qed
\end{cor}
\begin{cor}
Suppose that $n=2m+1$, $m\in \ZZ_{\ge 0}$. 
We have 
\begin{eqnarray*}
(\chi(E_i,E_j))_{i,j=a}^{\widetilde \mu_n-1+a}&=&(1-N)\cdot 
\frac{1-(N^{d_1})^{a_2}}{1-N^{d_1}}\cdot \frac{1-(N^{d_3})^{a_4}}{1-N^{d_3}}\cdot \dots \cdot \frac{1-(N^{d_{2m-1}})^{a_{2m}}}{1-N^{d_{2m-1}}}\\
&=&\prod_{i=0}^{n}\left(1-N^{d_i}\right)^{(-1)^{n-i+1}}=\frac{1}{\varphi_n(N)}=\chi_n.
\end{eqnarray*}
\qed
\end{cor}
In these corollaries, we used the fact that  $N^{d_n}=0$ due to $d_n> \widetilde \mu_n$.
\subsection{Grading shifts}
For each $a\in\ZZ$, denote by $\E_a$ the full triangulated subcategory $\langle E_a,\dots E_{\widetilde \mu_n-1+a}\rangle$ of ${\rm HMF}^{L_f}_{S}(f)$.
We shall show that $\E_a$ is closed under the grading shift $(\vec{l})$ for all $\vec{l}\in L_f$, which is done inductively. 
Since $L_f/\ZZ\vec{f}$ is generated by $\vec{x}_1$ and $(\vec{f})=T^2$, we only need to show that it is closed under $(\vec{x}_1)$.
\begin{lem}\label{lem:4.5}
Suppose that $n=2m$, $m\in \ZZ_{\ge 1}$. 
For each $i\in\ZZ$, set 
\[
E'_i:=\left(S/(x_1, x_2,x_4, \dots, x_{2m})\right)^{stab}(i\vec{x}_1).
\]
Then, there is an exact triangle 
\[
E_{i-1}\longrightarrow E_{i}\longrightarrow E'_i\longrightarrow TE_{i-1}.
\]
\end{lem}
\begin{pf}
This is straightforward from the definition of $E_i$ and Lemma~\ref{lem:4.1}.
\qed
\end{pf}
\begin{lem}\label{lem:4.6}
Suppose that $n=2m$, $m\in \ZZ_{\ge 1}$. 
Define $f'\in S':=\CC[x_2,\dots, x_n]$ by 
\[
f'=x_2^{a_2}x_3 + \cdots + x_{2m-1}^{a_{2m-1}}x_{2m} + x_{2m}^{a_{2m}}, 
\]
and set 
\[
\widetilde \mu'_{2m-1}:=\frac{\widetilde \mu_{2m}-1}{a_1}=a_2\dots a_{2m}-a_2\dots a_{2m-1}+\dots -a_2a_3+a_2-1.
\]
Then $(E'_{a_1+a},E'_{2a_1+a}\dots E'_{\widetilde \mu'_{2m-1}\cdot a_1+a})$ forms an exceptional collection for each $a\in\ZZ$.
If we further assume that Theorem~\ref{thm: main} holds for $n=2m-1$, 
then there is a triangulated equivalence 
\[
\langle E'_{a_1+a},E'_{2a_1+a},\dots, E'_{\widetilde \mu'_{2m-1}\cdot a_1+a} \rangle \cong {\rm HMF}^{L_{f'}}_{S'}(f').
\]
In particular, it is closed under the grading shift $(\vec{x}_2)$.
\end{lem}
\begin{pf}
Identify $L_{f'}$ with the subgroup of $L_f$ generated by $\vec{x}_2,\dots, \vec{x}_{2m}, \vec{f}$.
The following expression of $f$ for the description of the stabilization $E'_i$
\[
f=x_1(x_1^{a_1-1}x_2)+x_2(x_2^{a_2-1}x_3)+x_4(x_3^{a_3}+x_4^{a_4-1}x_5)+\dots +x_{2m}(x_{2m-1}^{a_{2m-1}}+x_{2m}^{a_{2m}-1})
\]
enables us to calculate the spaces of morphisms
\[
\bigoplus_{\vec{l}'\in L_{f'}}{\rm HMF}^{L_f}_{S}(f)(E'_i, E'_i(\vec{l}')),\quad \bigoplus_{\vec{l}'\in L_{f'}}{\rm HMF}^{L_f}_{S}(f)(E'_i, TE'_i(\vec{l}')),
\]
which turn out to be the same as those in Lemma~\ref{lem:4.2} for $f'$.
The rest is clear.
\qed
\end{pf}
\begin{cor}
Suppose that $n=2m$, $m\in \ZZ_{\ge 1}$ and that Theorem~\ref{thm: main} holds for $n=2m-1$. 
Then, $\E_a$ is closed under the grading shift $(\vec{x}_1)$.
\end{cor}
\begin{pf}
It is straightforward from Lemma~\ref{lem:4.5} that $E'_{a_1+a},E'_{2a_1+a},\dots, E'_{\widetilde \mu'_{2m-1}\cdot a_1+a}$ is 
in the subcategory since $a_1\ge 1$ and $\widetilde \mu'_{2m-1}\cdot a_1=\widetilde \mu_{2m}-1$.
Then, Lemma~\ref{lem:4.6} implies that $E'_{a}$ is also in the subcategory since $a_1\vec{x}_1=-\vec{x}_2$ in $L_f/\ZZ\vec{f}$.
Thus we see that $E_{a-1}\in \E_a$ by Lemma~\ref{lem:4.5} with $i=a$. 
Repeating this argument and using the fact that $d_{2m}\vec{x_1}=\vec{0}$ in $L_f/\ZZ\vec{f}$ and $(\vec{f})=T^2$, we obtain the statement.
\qed
\end{pf}
\begin{lem}\label{lem:4.8}
Suppose that $n=2m+1$, $m\in \ZZ_{\ge 1}$. 
For each $i\in\ZZ$ and $j=1,\dots, a_1$, set 
\[
E''_{i,j}:=\left(S/(x_1^j, x_3,x_5, \dots, x_{2m+1})\right)^{stab}(-i\vec{x}_1).
\]
Then, there is an exact triangle 
\[
E''_{i+1,j}\longrightarrow E''_{i,j+1}\oplus E''_{i+1,j-1}\longrightarrow E''_{i,j}\longrightarrow TE''_{i+1,j},
\]
where $E''_{i,0}:=0$ and $E''_{i,a_1+1}:=0$.
\end{lem}
\begin{pf}
It follows from a direct calculation.
\qed
\end{pf}
For each $i\in\ZZ$, set $E''_i:=E''_{i,a_1}$. Note that we have 
\[
E''_i\cong \left(S/(x_1^{a_1}+x_2^{a_2-1}x_3, x_3,x_5, \dots, x_{2m+1})\right)^{stab}(-i\vec{x}_1).
\]
\begin{lem}\label{lem:4.9}
Suppose that $n=2m+1$, $m\in \ZZ_{\ge 1}$. 
Define $f''\in S'':=\CC[x_3,\dots, x_n]$ by 
\[
f''=x_3^{a_3}x_4 + \cdots + x_{2m}^{a_{2m}}x_{2m+1} + x_{2m+1}^{a_{2m+1}}, 
\]
and set 
\[
\widetilde \mu''_{2m-1}:=\frac{\widetilde \mu_{2m+1}-a_1+1}{d_2}=a_3\dots a_{2m+1}-a_3\dots a_{2m}+\dots -a_3a_4+a_3-1.
\]
Then $(E''_{d_2+a_1-2+a},E''_{2d_2+a_1-2+a},\dots, E''_{\widetilde \mu''_{2m-1}\cdot d_2+a_1-2+a})$ forms an exceptional collection for each $a\in\ZZ$.
If we further assume that Theorem~\ref{thm: main} holds for $n=2m-1$, 
then there is a triangulated equivalence 
\[
\langle E''_{d_2+a_1-2+a},E''_{2d_2+a_1-2+a},\dots, E''_{\widetilde \mu''_{2m-1}\cdot d_2+a_1-2+a} \rangle \cong {\rm HMF}^{L_{f''}}_{S''}(f'').
\]
In particular, it is closed under the grading shift $(\vec{x}_3)$.
\end{lem}
\begin{pf}
Identify $L_{f''}$ with the subgroup of $L_f$ generated by $\vec{x}_3,\dots, \vec{x}_{2m+1}, \vec{f}$.
The following expression of $f$ for the description of the stabilization $E''_i$
\[
f=(x_1^{a_1}+x_2^{a_2-1}x_3)x_2 +x_3(x_3^{a_3-1}x_4)+x_5(x_4^{a_4}+x_5^{a_5-1}x_6)+\dots +x_{2m+1}(x_{2m}^{a_{2m}}+x_{2m+1}^{a_{2m+1}-1})
\]
enables us to calculate the spaces of morphisms 
\[
\bigoplus_{\vec{l}''\in L_{f''}}{\rm HMF}^{L_f}_{S}(f)(E''_i, E''_i(\vec{l}'')),\quad \bigoplus_{\vec{l}''\in L_{f'}}{\rm HMF}^{L_f}_{S}(f)(E''_i, TE''_i(\vec{l}'')),
\]
which turn out to be the same as those in Lemma~\ref{lem:4.2} for $f''$.
The rest is clear.
\qed\end{pf}
\begin{cor}
Suppose that $n=2m+1$, $m\in \ZZ_{\ge 1}$ and that Theorem~\ref{thm: main} holds for $n=2m-1$. 
Then, $\E_a$ is closed under the grading shift $(\vec{x}_1)$.
\end{cor}
\begin{pf}
Lemma~\ref{lem:4.8} implies that the objects $E''_{d_2+a_1-2+a}$, $E''_{2d_2+a_1-2+a}$, $\dots$, $E''_{\widetilde \mu''_{2m-1}\cdot d_2+a_1-2+a}$ are 
in the subcategory since $E''_{i,1}=E_i$ and $\widetilde \mu''_{2m-1}\cdot d_2+a_1-2+a=\widetilde \mu_{2m+1}-1+a$.
Then, it follows from Lemma~\ref{lem:4.9} that $E''_{a_1-2+a}$ is also in the subcategory since $d_2\vec{x}_1=\vec{x}_3$ in $L_f/\ZZ\vec{f}$.
Thus we see that $E_{a-1}\in \E_a$ by the use of Lemma~\ref{lem:4.8} with $i=a,\dots, a_1-2+a$. 
Repeating this argument and using the fact that $d_{2m}\vec{x_1}=\vec{0}$ in $L_f/\ZZ\vec{f}$ and $(\vec{f})=T^2$, we obtain the statement.
\qed
\end{pf}
\subsection{$\E_a$ contains $\CC^{stab}$}
Since $\E_a$ is closed under the grading shift $(\vec{l})$ for all $\vec{l}\in L_f$, it is enough to show that 
$\CC^{stab}(\vec{l})\in\E_a$ for some $\vec{l}\in L_f$.
If $n=1$, then it is obvious that $\E_a$ contains $\CC^{stab}(-a\vec{x}_1)$.
Suppose that $n=2m+1$, $m\in \ZZ_{\ge 1}$. 
For each $k=1,\dots, m$, there is an exact triangle 
\begin{eqnarray*}
& &\left(S/(x_1, x_2,\dots, x_{2k-2}, x_{2k-1}, x_{2k+1},x_{2k+3},\dots, x_{2m+1})\right)^{stab}(-\vec{x}_{2k})\\
& &\longrightarrow \left(S/(x_1, x_2,\dots, x_{2k-2}, x_{2k-1}, x_{2k+1},x_{2k+3},\dots, x_{2m+1})\right)^{stab}\\
& &\longrightarrow \left(S/(x_1, x_2,\dots, x_{2k-2}, x_{2k-1}, x_{2k}, x_{2k+1},x_{2k+3},\dots, x_{2m+1})\right)^{stab}\\ 
& &\longrightarrow T\left(S/(x_1, x_2,\dots, x_{2k-2}, x_{2k-1}, x_{2k+1},x_{2k+3},\dots, x_{2m+1})\right)^{stab}(-\vec{x}_{2k}).
\end{eqnarray*}
Note that $-\vec{x}_2-\vec{x}_4-\dots -\vec{x}_{2m}=(d_1+d_3+\dots +d_{2m-1})\vec{x}_1$ in $L_f/\ZZ\vec{f}$. 
Since $a_i\ge 2$ for all $i=1,\dots, 2m+1$, it follows that 
\begin{eqnarray*}
\widetilde \mu_{2m+1}&=&d_{2m-1} a_{2m}(a_{2m+1}-1)+\widetilde \mu_{2m-1}\\
&>&d_{2m-1}+\widetilde \mu_{2m-1}\\
&>&d_{2m-1}+d_{2m-1}+\widetilde \mu_{2m-3}\\
& >& \dots > d_{2m-1}+d_{2m-1}+\dots +d_3+d_1.
\end{eqnarray*}
Hence, the above exact triangle yields that $\CC^{stab}(\vec{l})\in \E_a$ for some $\vec{l}\in L_f$.
Similarly, suppose that $n=2m$, $m\in \ZZ_{\ge 1}$. 
For each $k=1,\dots, m$, there is an exact triangle 
\begin{eqnarray*}
& &\left(S/(x_1, x_2,\dots, x_{2k-3}, x_{2k-2}, x_{2k}, x_{2k+2},\dots, x_{2m})\right)^{stab}(-\vec{x}_{2k-1})\\
& &\longrightarrow \left(S/(x_1, x_2,\dots, x_{2k-3}, x_{2k-2}, x_{2k}, x_{2k+2},\dots, x_{2m})\right)^{stab}\\
& &\longrightarrow \left(S/(x_1, x_2,\dots, x_{2k-3}, x_{2k-2}, x_{2k-1}, x_{2k}, x_{2k+2},\dots, x_{2m})\right)^{stab}\\ 
& &\longrightarrow T\left(S/(x_1, x_2,\dots, x_{2k-3}, x_{2k-2}, x_{2k}, x_{2k+2},\dots, x_{2m})\right)^{stab}(-\vec{x}_{2k-1}).
\end{eqnarray*}
Note that $-\vec{x}_1-\vec{x}_3-\dots -\vec{x}_{2m-1}=-(d_0+d_2+\dots +d_{2m-2})\vec{x}_1$ in $L_f/\ZZ\vec{f}$. 
Since $a_i\ge 2$ for all $i=1,\dots, 2m+1$, it follows that 
\begin{eqnarray*}
-\widetilde \mu_{2m}&=&-d_{2m-2} a_{2m-1}(a_{2m}-1)-\widetilde \mu_{2m-2}\\
&<&-d_{2m-2}-\widetilde \mu_{2m-2}\\
&<&-d_{2m-2}-d_{2m-4}-\widetilde \mu_{2m-4}\\
& <& \dots < -d_{2m-2}-d_{2m-3}-\dots -d_2-d_0.
\end{eqnarray*}
Hence, the above exact triangle yields that $\CC^{stab}(\vec{l})\in \E_a$ for some $\vec{l}\in L_f$.
Thus we have finished the proof of Theorem~\ref{thm: main}.

\end{document}